\documentclass[a4paper,oneside,12pt]{article}
\usepackage[english, russian]{babel}
\usepackage{amsmath}
\usepackage{amsfonts}
\usepackage{amsthm}

\theoremstyle{plain} \newtheorem{theorem}{\qquad Теорема}
\theoremstyle{remark} \newtheorem*{remark}{\qquad Замечание}
\numberwithin{equation}{section}
\everymath{\displaystyle}

\begin{document}
УДК 517.955.8

\begin{center}
П.В. Бабич, В.Б. Левенштам, С.П. Прика

\textbf{Восстановление быстро осциллирующего источника в уравнении теплопроводности по асимптотике решения}
\end{center}

\section*{Введение}
\label{intro}

\qquad В работе рассматриваются некоторые асимптотические задачи для уравнения теплопроводности в прямоугольнике с быстро осциллирующим по времени источником вида $f(x,t) r(t,\omega t)$ ($r$ -- $2\pi$-периодична по второй переменной, $\omega \gg 1$) и однородными начально-краевыми условиями. В каждой задаче источник не известен (точнее, не известен один из сомножителей $f,r$ или оба, но известны классы функций, которым они принадлежат) и требуется восстановить его по тем или иным сведениям об асимптотике решения. Исследуются четыре случая:
1) неизвестной функции $r(t,\tau)$;
2) неизвестной функции $f(x,t) \equiv f(x)$;
3) неизвестных $f(x)$ и $r_1 (t,\tau) = r(t,\tau) - r_0 (t)$, где $r_0 (t) = (2\pi)^{-1} \int\limits_0^{2\pi} r(t,\tau) d\tau$
и, наконец, случай 4) неизвестных $f(x)$ и $r(t,\tau)$. В случае 4) при этом предполагается, что $f(x)$ разлагается в сумму конечного числа $N$ гармоник.
Для задач 1)-3) ставятся соответственно следующие дополнительные условия: 1) задана двучленная асимптотика (при $\omega\to\infty$) решения, вычисленная в некоторой точке $x=x_0$, 2) задан главный член асимптотики решения в некоторой точке $t = t_0$, 3) одновременно заданы условия 1), 2) в точках $x = x_0, t = t_0$. В задаче 4) заданы условия 3) и, кроме того, известен главный член асимптотики решения $u_0 (x,t)$ для $N-1$ различных точек $x = x_1, x_2, ... , x_{N-1}$ в некоторой окрестности точки $t=t_0$. Другими словами, заданы функции $\alpha_j (t), t \in (t_0-\delta, t_0 + \delta) \equiv I_{\delta}, \delta >0, j = \overline{1,N}$, такие что $\alpha_j (t) = u_0 (x_j, t), t \in I_{\delta}$.

Данная статья стимулирована важными работами А.~М.~Денисова \cite{babich-bib1, babich-bib2}, в которых исследован целый ряд различных обратных коэффициентных задач для параболических и гиперболических уравнений (без высокочастотных осцилляций). В частности, в \cite{babich-bib1} рассмотрена обратная задача типа 1) для уравнения теплопроводности с источником вида $f(x)r(t)$, где функция $r$ не известна (в \cite{babich-bib1} это первая задача в соответствующей паре: параболическое+гиперболическое уравнения); дополнительно задано решение начально-краевой задачи, вычисленное в точке $x=x_0$. Для задач в нашей ситуации, т.е. соответствующих уравнениям с быстро осциллирующими данными, постановка дополнительного условия на все решение (пусть даже вычисленное в одной точке $x = x_0$) неестественно (слишком обременительно). Естественно ставить такое условие лишь на первые члены асимптотики. Сколько этих членов должно быть задействовано для определения неизвестных следует выяснять при анализе асимптотического разложения решения исходной задачи. Таким путем были получены указанные выше дополнительные условия в задачах 1)-4). Следует сказать, что очень близкая к результату 1) теорема установлена нами в \cite{babich-bib3}, где задача поставлена немного иначе и $f(x,t)$ не зависит от $t$.

Отметим, что задачи с быстро осциллирующими по времени данными моделируют многие физические (и иные) процессы (в частности, связанные с высокочастотными механическими воздействиями на среду). К ним относятся, например: диффузия вещества в среде, подверженной высокочастотным вибрациям; коротковолновая дифракция акустических и продольных волн; конвекция жидкости в поле быстро осциллирующих сил \cite{babich-bib31}--\cite{babich-bib5} и др.

Заметим еще, что все четыре рассматриваемые в данной статье задачи, относящиеся к одномерному уравнению теплопроводности, решены в настоящее время одним из соавторов и для многомерного случая. При этом широко и эффективно используемые в одномерной ситуации ряды Фурье в многомерном случае используются минимально. Объясняется это тем, что обоснование метода Фурье в многомерном случае традиционно связано с требованиями чрезмерно завышенной гладкости данных задачи. Указанные многомерные результаты предполагается опубликовать в отдельной статье.

Отметим, наконец, что в настоящей работе мы налагаем на данные исходной задачи условия, обеспечивающие достаточную гладкость ее решения в замкнутом прямоугольнике $\overline{\Pi}$. Можно решить рассматриваемые здесь задачи (в одномерной и многомерной ситуациях) и при менее жестких требованиях -- когда решение содержит так называемый пространственный погранслой \cite{babich-bib4}. Вопросам построения и обоснования асимптотик решений с погранслоем начально-краевых задач и задач о периодических по времени решениях для различных параболических уравнений посвящены работы \cite{babich-bib6}-\cite{babich-bib7} др. Соответствующие задачи о восстановлении неизвестных исходных данных мы планируем исследовать позже.

\section{Неизвестный сомножитель источника зависит только от времени}
\label{sect1}

\subsection*{$1.1^{\circ}$ Двучленная асимптотика}

\qquad Символом $\Pi$ обозначим открытый прямоугольник ${\{ (x,t): 0 < x < \pi; 0 < t < T \}}, T = const > 0 $, а символом  $ \Gamma $ - часть его границы, состоящую из нижнего основания и боковых сторон. Рассмотрим параболическую начально-краевую задачу с большим параметром $\omega$:
\begin{equation}
\frac{\partial u}{\partial t} = \frac{\partial^2 u}{\partial x^2} + f(x,t) r(t, \omega t), (x,t) \in \Pi,
\label{s1}
\end{equation}
\begin{equation}
\left. u \right|_{\Gamma} = 0.
\label{s2}
\end{equation}
Относительно функций $f(x,t)$ и $r(t,\tau)$, определенных и непрерывных соответственно на множествах $(x,t) \in [0,\pi] \times [0,T] \equiv S_0$ и $(t,\tau) \in [0,T] \times [0,\infty) \equiv S_1$, сделаем следующие предположения. Функция $r(t,\tau)$ -- 2$\pi$-периодична по $\tau$; представим ее в виде $r(t,\tau) = r_0(t) + r_1(t,\tau)$, где $r_1(t,\tau)$ имеет нулевое среднее по второй переменной:
\begin{equation*}
\left\langle r_1(t,\cdot) \right\rangle = \left\langle r_1(t,\tau) \right\rangle_{\tau} \equiv \frac{1}{2\pi} \int\limits_0^{2\pi} r_1(t,\tau) d\tau = 0.
\label{r_1}
\end{equation*}
Пусть $f \in C^{3,0}(S_0), r_0 \in C^0 ([0,T]), r_1 \in C^{1+\alpha,0} (S_1), \alpha\in (0,1),$
\begin{equation}
f(0,t) = f(\pi,t) = f''(0,t) = f''(\pi,t) = 0.
\label{f0}
\end{equation}
Здесь и ниже символом $C^{l,m}(S)$, где $l,m$ -- неотрицательные числа, мы обозначаем обычные гельдеровы пространства функций. Функцию $r(t,\tau)$, удовлетворяющую указанным выше условиям, будем, для краткости, называть функцией класса (А). В этом пункте будет построена и обоснована двучленная асимптотика в $C (\overline{\Pi})$ (относительно малого параметра $\omega^{-1}$) решения задачи (\ref{s1})-(\ref{s2}).

Решение задачи (\ref{s1})-(\ref{s2}) запишем в виде:
\begin{equation}
u_{\omega} (x,t) = U_{\omega} (x,t) + W_{\omega} (x,t), \omega \gg 1,
\label{uw}
\end{equation}
где
\begin{equation}
U_{\omega} (x,t) = u_0 (x,t) + \omega^{-1} \bigl[ u_1 (x,t) + v_1 (x,t,\omega t)\bigr], \omega \gg 1,
\label{uww}
\end{equation}
\begin{equation}
u_0 (x,t) =  \sum_{n=1}^{\infty} \sin{nx} \int\limits_0^t e^{-n^2(t-s)} f_n (s) r_0 (s) ds, \; f_n (t) = \frac{2}{\pi} \int\limits_0^{\pi} f(s,t) \sin{ns} ds,
\label{u_0}
\end{equation}
\begin{equation}
v_1 (x,t,\tau) = f(x,t) \Biggl[ \int\limits_0^{\tau} r_1 (t,s) ds - \left\langle \int\limits_0^{\tau} r_1 (t,s) ds \right\rangle_{\tau} \Biggr],
\label{v_1}
\end{equation}
\begin{equation}
u_1 (x,t) = \left\langle \int\limits_0^{\tau} r_1 (0,s) ds \right\rangle_{\tau} \cdot \sum_{n=1}^{\infty} f_n (0) \sin{nx} e^{-n^2 t}.
\label{u_1}
\end{equation}
Выражения (\ref{u_0})-(\ref{u_1}) будут выведены при доказательстве следующей теоремы.

\begin{theorem}
Решение $u_{\omega}(x,t)$ задачи (\ref{s1})-(\ref{s2}) представимо в виде (\ref{uw})-(\ref{u_1}), где
\begin{equation}
\bigl\| W_{\omega} (x,t) \bigr\|_{C(\overline{\Pi})} = o(\omega^{-1}), \omega \to \infty. \label{th1}
\end{equation} \label{theor1}
\end{theorem}

\begin{remark}
Можно построить и обосновать полную ($N$-членную при любом фиксированном $N$) асимптотику решения задачи \eqref{s1}-\eqref{s2} по норме $C(\overline{\Pi})$ для бесконечно дифференцируемых (достаточно гладких) функции $r(t,\tau)$ и $f(x)$, удовлетворяющих условиям типа \eqref{f0} любого (достаточно высокого) порядка. Мы этого здесь не делаем, поскольку для решения задачи о восстановлении источника нам понадобится лишь двучленная асимптотика.
\end{remark}

\subsection*{$1.2^{\circ}$ Восстановление источника}
\label{sect3}

\qquad Пусть в начально-краевой задаче \eqref{s1}-\eqref{s2} фигурирует та же, что и в п.$1^{\circ}$, функция $f(x,t)$, а функция $r(t,\tau)$ класса (А) не известна. Пусть $x_0 \in (0,\pi)$ -- точка, в которой $f(x_0,t) \neq 0$. Введем в рассмотрение функции $\varphi_0 (t), \varphi_1(t), \varphi_2 (t,\tau)$, удовлетворяющие следующим условиям:
\begin{equation}
\varphi_0 \in C^1([0,T]), \; \varphi_0 (0) = 0,
\label{eq19'}
\end{equation}
$\varphi_2 \in C^{1+\alpha,1} (S_2)$ $2\pi$-периодична по $\tau$ с нулевым средним: $\left\langle \varphi_2 (t,\cdot) \right\rangle = 0,$
\begin{equation*}
\varphi_1 (t) = \frac{1}{f(x_0,0)} \left\langle \int\limits_0^{\tau} \frac{\partial}{\partial \tau} \varphi_2 (0,s) ds \right\rangle_{\tau} \cdot \sum_{n=1}^{\infty} f_n (0) \sin{nx} e^{-n^2 t}.
\end{equation*}

Наша задача -- назовем ее задачей 1 -- состоит в определении такой функции $r (t,\tau)$ класса (А), для которой решение $u_{\omega} (x,t)$ задачи \eqref{s1}-\eqref{s2} удовлетворяет условию
\begin{equation}
\left\| u_{\omega} (x_0,t) - \left[ \varphi_0 (t) + \frac{1}{\omega}\bigl( \varphi_1 (t) + \varphi_2 (t,\omega t) \bigr) \right] \right\|_{C([0,T])} = o (\omega^{-1}), \; \omega \to \infty.
\label{babich-inv}
\end{equation}

\begin{theorem}
Для любой тройки функций $\varphi_0, \varphi_1$ и $\varphi_2$ и точки $x_0$, удовлетворяющих указанным выше условиям, существует единственная функция $r$ класса (А), для которой решение $u_{\omega} (x,t)$ задачи \eqref{s1}-\eqref{s2} удовлетворяет условию \eqref{babich-inv}. Нахождение функции $r$ сводится к решению уравнения Вольтерры второго рода.
\label{babich-t4}
\end{theorem}

\subsection*{$1.3^{\circ}$ Доказательство основных результатов параграфа}
\label{sect2}

\noindent\textbf{Доказательство теоремы \ref{theor1}.}

Подставим выражение $u_{\omega}$ (\ref{uw}) в равенства (\ref{s1})-(\ref{s2}):
\begin{equation}
\left\{
\begin{array}{l}
     \frac{\partial U_{\omega}}{\partial t} + \frac{\partial W_{\omega}}{\partial t} = \frac{\partial^2 U_{\omega}}{\partial x^2} + \frac{\partial^2 W_{\omega}}{\partial x^2} + f(x)r(t,\omega t), \omega \gg 1,\\
     \left. \bigl( U_{\omega} + W_{\omega} \bigr) \right|_{\Gamma} = 0.
\end{array}
\right.
\label{uuw}
\end{equation}
Учитывая (\ref{uww}), получим
\begin{equation}
\left\{
\begin{array}{l}
    \displaystyle{\frac{\partial u_0}{\partial t} + \frac{\partial v_1}{\partial \tau} + \omega^{-1} \left[ \frac{\partial u_1}{\partial t} + \frac{\partial v_1}{\partial t} \right] + \frac{\partial W_{\omega}}{\partial t} =}\\
       \qquad \qquad \displaystyle{ = \frac{\partial^2 u_0}{\partial x^2} + \omega^{-1} \left[\frac{\partial^2 u_1}{\partial x^2} + \frac{\partial^2 v_1}{\partial x^2} \right] + \frac{\partial^2 W_{\omega} }{\partial x^2} + f(x)r(t,\tau), } \\
    \displaystyle{\left. \left( u_0 + \omega^{-1} (u_1 + v_1) + W_{\omega} \right) \right|_{\Gamma} = 0, \; \tau = \omega t, \; \omega \gg 1.}
\end{array}
\right.
\label{suw}
\end{equation}
Приравняем в последних равенствах коэффициенты при степенях $\omega^{-k}, \; (k = 0,1)$, а затем применим к полученным уравнениям операцию усреднения $ \langle ...\rangle_{\tau} $ по $\tau = \omega t$. В результате придем к следующим задачам:
\begin{equation}
        \left\{
        \begin{array}{l}
           \displaystyle{\frac{\partial u_0}{\partial t} = \frac{\partial^2 u_0}{\partial x^2} + f(x,t) r_0 (t),} \\
           \displaystyle{\left. u_0 \right|_{\Gamma} = 0} ,
        \end{array}
        \right.
        \label{sA}
\end{equation}
  \begin{equation}
        \left\{
        \begin{array}{l}
           \displaystyle{\frac{\partial v_1}{\partial \tau} = f(x,t) r_1 (t,\tau), }\\
           \displaystyle{\left\langle v_1 (x,t,\tau) \right\rangle_{\tau} = 0,}
        \end{array}
        \right.
        \label{sV1}
  \end{equation}
\begin{equation}
       \left\{
       \begin{array}{l}
            \displaystyle{\frac{\partial u_1}{\partial t} = \frac{\partial^2 u_1}{\partial x^2}, } \\
            \displaystyle{\left. u_1 \right|_{t = 0} = \left. - v_1 \right|_{\begin{smallmatrix} t = 0\\\tau = 0 \end{smallmatrix}}, }\\
            \displaystyle{\left. u_1 \right|_{x = 0} = \left. u_1 \right|_{x = \pi} = 0, }
       \end{array}
       \right.
       \label{sU1}
\end{equation}
Отметим, что в системе (\ref{sU1}) учтены очевидные соотношения: \linebreak ${\left. v_1 \right|_{x=0,\pi} = 0}$. Кроме того, отметим, что в силу (\ref{sA})-(\ref{sU1}) функции $u_0, v_1$ и $u_1$ имеют вид (\ref{u_0}), (\ref{v_1}) и (\ref{u_1}) соответственно. Это устанавливается с помощью элементарного применения рядов Фурье.

Из равенств (\ref{suw})-(\ref{sU1}) находим
\begin{equation}
\left\{
\begin{array}{l}
   \displaystyle \frac{\partial W_{\omega}}{\partial t} = \frac{\partial^2 W_{\omega}}{\partial x^2} + \omega^{-1} \left[ \frac{\partial^2 v_1}{\partial x^2} - \frac{\partial v_1}{\partial t}\right], \\
   \displaystyle{\left. W_{\omega} \right|_{\Gamma} = 0.}
\end{array}
\right.
\label{esw}
\end{equation}
Решение задачи \eqref{esw} в силу \eqref{v_1} имеет вид
\begin{multline}
W_{\omega} (x,t) = \\
  \frac{1}{\omega} \sum_{n = 1}^{\infty} \sin nx \left[ f''_n \int\limits_0^t e^{-n^2 (t-s)}  p(s,\omega s) ds - f_n \int\limits_0^t e^{-n^2 (t-s)} \left. \frac{\partial}{\partial t} p(t,\tau) \right|_{\begin{smallmatrix} t = s\\\tau = \omega s \end{smallmatrix}} ds \right] \equiv\\
    W_{\omega,1}(x,t) + W_{\omega,2}(x,t),
\label{solW}
\end{multline}
где
\begin{equation}
f''_n = \frac{1}{\pi} \int\limits_0^{\pi} f''(s) \sin ns ds, \; p(t,\tau) = \int\limits_0^{\tau} r_2 (t,s) ds - \left\langle \int\limits_0^{\tau} r_2 (t,s) ds \right\rangle_{\tau}.
\label{wf}
\end{equation}

Для доказательства теоремы 1 достаточно установить соотношения:
\begin{equation}
\| W_{\omega,1} \|_{C(\overline{\Pi})} = o(\omega^{-1}), \; \| W_{\omega,2} \|_{C(\overline{\Pi})} = o(\omega^{-1}), \; \omega \to \infty,
\label{W12}
\end{equation}
Поскольку рассуждения при их выводе совершенно аналогичны, то ограничимся выводом первого из них. В силу того, что $f \in C^3 ([0,\pi])$ и $f^{(2k)} (0) = f^{(2k)} (\pi) = 0, k=0,1,$ функцию $f''(x)$ можно продолжить нечетным, $2\pi$-периодическим образом с сохранением класса гладкости $C^3$ на всю прямую $x \in \mathbb{R}$. Оставив за продолжением прежнее обозначение, получим
\begin{equation*}
f''(x) = \sum_{n = 1}^{\infty} f''_n \sin nx, x \in \mathbb{R}.
\end{equation*}
Поскольку $f'' \in C^1 (\mathbb{R})$, то
\begin{equation*}
f''_n = \frac{a_n}{n}, \; \sum_{n = 1}^{\infty} a_n^2 < \infty.
\end{equation*}
Из этих соображений и неравенства Коши-Буняковского следует, что нам теперь достаточно доказать равномерную относительно $n \in \mathbb{N}$ асимптотическую оценку
\begin{equation}
\left\| \int\limits_0^t e^{-n^2 (t-s)}  p(s,\omega s) ds \right\|_{C([0,1])} = o (1), \omega \to \infty.
\label{op}
\end{equation}
Проведем ее в три этапа. Пусть $\varepsilon$ -- произвольное положительное число. На первом этапе подберем $\delta > 0$ столь малым, что при всех $t \in [0,1], n \in \mathbb{N}$ и $\omega > 0$
\begin{equation}
\left| \int\limits_{t_0}^t e^{-n^2 (t-s)} p(s,\omega s) ds \right| < \frac{\varepsilon}{3}, t_0 = \max (0,t-\delta).
\label{estW1}
\end{equation}
На втором этапе при $t_0 = t-\delta > 0$ подберем натуральное число $n_0$ столь большим, что при $n \geq n_0, t \geq \delta$ и $\omega > 0$
\begin{equation}
\left| \int\limits_0^{t-\delta} e^{-n^2 (t-s)} p(s,\omega s) ds \right| < \frac{\varepsilon}{3}.
\label{estW2}
\end{equation}
На третьем этапе разобьем участок $[0,t-\delta], t > \delta,$ на $m$ равных частей $[t_j, t_{j+1}), j = 0,1, \ldots, m-1$. Воспользуемся равенством
\begin{multline*}
\int\limits_0^{t-\delta} e^{-n^2 (t-s)}  p(s,\omega s) ds = \\
 \sum_{j = 0}^{m-1} \left[ \int\limits_{t_j}^{t_{j+1}} e^{-n^2 (t-s)} p(s,\omega s) ds - \int\limits_{t_j}^{t_{j+1}} e^{-n^2 (t - t_j)} p (t_j,\omega s) ds \right] + \\
   \sum_{j = 0}^{m-1} \int\limits_{t_j}^{t_{j+1}} e^{-n^2 (t - t_j)} p (t_j,\omega s) ds = S_1 + S_2.
\end{multline*}
Выберем затем $m$ столь большим, что при всех натуральных $n \leq n_0$ и $\omega > 0$
\begin{equation}
| S_1 | < \frac{\varepsilon}{6}.
\end{equation}
Далее, в силу равенства
\begin{equation*}
\sum_{j = 0}^{m-1} \int\limits_{t_j}^{t_{j+1}} e^{-n^2 (t - t_j)} p (t_j,\omega s) ds = \frac{1}{\omega} e^{-n^2 (t-t_j)} \left[ \int\limits_{0}^{\omega t_{j+1}} p(t_j,\tau) d\tau - \int\limits_{0}^{\omega t_j} p(t_j,\tau) d\tau \right]
\end{equation*}
и того факта, что функция $p(t,\tau)$, имеет нулевое среднее по второй переменной, найдется такое $\omega_0 > 0$, что для любого $\omega > \omega_0$
\begin{equation}
| S_2 | < \frac{\varepsilon}{6}.
\label{ests2}
\end{equation}
Из соотношений \eqref{estW1} -- \eqref{ests2} следует \eqref{op}. Теорема 1 доказана.

\noindent\textbf{Доказательство теоремы 2.}

Из теоремы \ref{theor1} следует, что при заданной функции $r(t,\tau)$ класса (А) задача (\ref{s1})-(\ref{s2}) имеет единственное классическое решение, представимое в виде:
\begin{equation*}
u_{\omega} (x,t) = u_0 (x,t) + \omega^{-1} \bigl[ u_1 (x,t) + v_1 (x,t,\omega t)\bigr] + W_{\omega} (x,t),
\end{equation*}
\begin{equation*}
\| W_{\omega} \|_{C (\overline{\Pi})} = o(\omega^{-1}), \; \omega \to \infty.
\end{equation*}
Пусть функция $r(t,\tau)$ является решением задачи 1 и $u_{\omega}$ -- отвечающее ему решение задачи \eqref{s1}-\eqref{s2}. В силу \eqref{th1},(\ref{babich-inv}) равномерно по $t \in [0,1]$
\begin{multline}
u_0 (x_0,t) + \omega^{-1} \bigl[ u_1 (x_0,t) + v_1 (x_0,t,\omega t)\bigr] = \\
   \varphi_0 (t) + \omega^{-1} \bigl[ \varphi_1 (t) + \varphi_2 (t, \omega t) \bigr] + o (\omega^{-1}), \; \omega \gg 1.
\label{inv-prove}
\end{multline}
%Продифференцировав (\ref{inv-prove}) по $t$, получим эквивалентное (\ref{inv-prove}) уравнение
%\begin{multline}
%\frac{\partial}{\partial t}u_0(x_0,t) + \left.\frac{\partial}{\partial \tau}v_1(x_0,t,\tau)\right|_{\tau = \omega t} + \omega^{-1} \left[ \frac{\partial}{\partial t}u_1(x_0,t) + \left. \frac{\partial}{\partial t}v_1(x_0,t,\tau)\right|_{\tau = \omega t} \right] + \frac{\partial}{\partial t}W_{\omega} (x_0,t) = \\
%     \varphi_0' (t) + \left.\frac{\partial}{\partial \tau} \varphi_1 (t,\tau)\right|_{\tau = \omega t} + \omega^{-1} \left.\frac{\partial}{\partial t} \varphi_1 (t,\tau)\right|_{\tau = \omega t} + \frac{d}{d t} \nu_{\omega} (t).
%\label{eq7}
%\end{multline}
Приравнивая в \eqref{inv-prove} коэффициенты при одинаковых степенях $\omega$, проводя затем усреднение по $\tau$ и дифференцируя полученные равенства, найдем
\begin{eqnarray}
\frac{\partial}{\partial t}u_0(x_0,t) = \varphi_0' (t), \label{eq80} \\
\frac{\partial}{\partial \tau}v_1(x_0,t,\tau) = \frac{\partial}{\partial \tau} \varphi_1 (t,\tau). \label{eq90}
\end{eqnarray}
В силу п. $1^{\circ}$ $\frac{\partial v_1 (x,t)}{\partial \tau} = f(x,t)r_1(t,\tau)$. Подставляя это выражение в (\ref{eq90}), получим
\begin{equation}
f(x_0,t)r_1(t,\tau) = \frac{\partial}{\partial \tau} \varphi_1 (t,\tau).
\label{eq9}
\end{equation}
Функция $u_0(x,t)$, согласно п. $1^{\circ}$, определяется равенством (\ref{u_0}), продифференцировав по $t$ которое, найдем
\begin{equation}
\frac{\partial}{\partial t}u_0(x_0,t) = f(x_0,t)r_0(t) + \int_0^t K(t,s) r_0(s) \, ds,
\label{eq81}
\end{equation}
где функция
\begin{equation}
K(t,s) = - \sum\limits_{n=1}^{\infty} n^2 f_n (s) e^{-n^2 (t-s)} \sin n x_0.
\label{kt}
\end{equation}

В силу условий, наложенных в п. \ref{sect1} на функцию $f$, $K(t,s)$ непрерывна. Из (\ref{eq80}), (\ref{eq81}) следует уравнение Вольтерра II
\begin{equation}
f(x_0,t)r_0(t) + \int_0^t K(t,s) r_0(s) \, ds = \varphi_0' (t),
\label{volt}
\end{equation}
которое имеет единственное непрерывное решение $r_0(t)$. Из уравнения (\ref{eq9}) функция $r_1$ также определяется единственным образом:
\begin{equation}
r_1 (t,\tau) = \frac{1}{f(x_0,t)} \frac{\partial}{\partial \tau} \varphi_2 (t,\tau),
\label{r1}
\end{equation}
и в силу условий, наложенных на функцию $\varphi_1$, принадлежит пространству $C^{1+\alpha,0}(S)$.

Так как найденная функция $r(t,\tau) = r_0(t) + r_1 (t,\tau)$ удовлетворяет условиям п. $1^{\circ}$, то для нее имеет место теорема 1, и, в частности, решение задачи \eqref{s1}-\eqref{s2}, представимо в виде \eqref{uw},\eqref{uww}. Покажем, что для функции $u_{\omega} (x,t)$ будет выполнено условие \eqref{babich-inv}. Для этого достаточно установить равенства:
\begin{equation*}
u_0 (x_0,t) = \varphi_0 (t), \; v_1 (x_0,t,\tau) = \varphi_2 (t,\tau), \; u_1 (x_0,t) = \varphi_1 (t).
\end{equation*}
Из предыдущей части доказательства теоремы 2 нам известно, что функция $r_0 (t)$ является решением уравнения \eqref{volt}. Из \eqref{eq81},\eqref{volt} следует равенство $\frac{\partial}{\partial t} u_0 (x_0,t) = \varphi'_0 (t)$. Учитывая, что $u_0 (x_0,0) = \varphi_0 (0) = 0$, получим $u_0 (x_0,t) = \varphi_0 (t)$. Функция $r_1$ определяется равенством \eqref{eq9}. Подставляя ее в \eqref{sV1} при $x=x_0$ и учитывая, что функции $\varphi_2 (t, \tau)$ и $v_1 (x,t,\tau)$, являются $2\pi$-периодическими с нулевым средним по $\tau$, получим $v_1 (x_0,t,\tau) = \varphi_2 (t,\tau)$. Отсюда, в силу \eqref{u_1}, находим
\begin{equation*}
u_1 (x_0,t) = \frac{1}{f(x_0)} \left\langle \int\limits_0^{\tau} \frac{\partial}{\partial \tau} \varphi_2 (0,s) ds \right\rangle_{\tau} \cdot \sum_{n=1}^{\infty} f_n  \sin{nx} e^{-n^2 t} = \varphi_1 (t).
\end{equation*}
Теорема доказана.

\section{Неизвестный сомножитель источника зависит от пространственной переменной}
\label{sect2}

\subsection*{$2.1^{\circ}$ Главный член асимптотики}

\qquad Пусть $\Pi$ и $\Gamma$ - те же, что в $\S 1$. Рассмотрим задачу \eqref{s1}-\eqref{s2} с $f(x,t) \equiv f(x)$. Предположим, что функция $f \in C^1([0,\pi]), f(0) = f(\pi) = 0$, а $r(t,\tau) = r_0 (t) + r_1 (t,\tau)$, где $r_0 \in C ([0,T]), r_1 \in C^{\alpha,0} (S), \alpha \in (0,1),$ и $r_1 (t,\tau)$ -- 2$\pi$-периодична по $\tau$ с нулевым средним
\begin{equation}
u_0 (x,t) =  \sum_{n=1}^{\infty} f_n \sin{nx} \int\limits_0^t e^{-n^2(t-s)} r_0 (s) ds \equiv \sum_{n = 1}^{\infty} f_n \Lambda_n \sin{nx}.
\label{eq2.11}
\end{equation}
Из соотношений: $f_n = n^{-1}b_n$, где $\sum_{n=1}^{\infty}b_n^2 < \infty,$ и неравенства Коши-Буняковского следует, что $u_0(x,t)$ - классическое решение задачи \eqref{s1}-\eqref{s2} при $r = r_0.$

\begin{theorem}
Справедлива асимптотическая формула
\begin{center}
$\bigl\|u_{\omega} - u_0 \bigr\|_{C(\overline{\Pi})} = o(1), \; \omega \to \infty.$
\end{center}
где $u_{\omega}$ - решение задачи (2.1), (2.2).
\end{theorem}

\subsection*{$2.2^{\circ}$ Восстановление источника}

\qquad При постановке рассматриваемой в данном пункте задачи 2 будем учитывать следующий простой факт. Если функция $r_0 \in C^{1},$ а функция $f \in C^{2}([0, \pi])$ и $f(0) = f(\pi)$, то для любого $t_0 \in (0,1]$ функция $u_0 (x,t_0) \equiv \psi (x)$ удовлетворяет условиям:
\begin{equation}
\psi (x) \in C^{4}([0, \pi]), \; \psi^{(2k)}(0) = \psi^{(2k)}(\pi) = 0, k = 0,1,2.
\label{eq2.1}
\end{equation}
Этот факт вытекает из следующих справедливых при любом фиксированном $t > 0$ простых соотношений:
\begin{equation*}
\Lambda'_n (t) = e^{- n^2 t} r_0 (0) - \int\limits_0^t e^{-n^2(t-s)} r'(s) ds = O(n^{-2}), n \to \infty,
\end{equation*}
\begin{equation*}
\frac{\partial^4 u_0}{\partial x^4} = \frac{\partial^2}{\partial x^2} \left( \frac{\partial u_0}{\partial t} - f(x) r_0 (t) \right) = \frac{\partial^2}{\partial x^2} \sum_{n = 1}^{\infty} f_n \Lambda'_n (t) \sin{x} -  f''(x) r_0 (t).
\end{equation*}

Итак, рассмотрим задачу \eqref{s1}-\eqref{s2} с $f(x,t) \equiv f(x)$. Будем считать, что функция $r(t, \tau)$ известна и удовлетворяет условиям $\S 1$, а функция $f \in C^1([0,\pi])$ не известна. Введем в рассмотрение точку $t_0 \in (0,T]$, в которой $r_0(t_0) \neq 0$. При этом естественном предположении справедливо следующее утверждение о числах $\Lambda_n \equiv \Lambda_n(t_0), n = 1,2,\ldots$ (см. \eqref{eq2.11}).

\textbf{Лемма.} \emph{Среди чисел $\Lambda_n$ может быть лишь конечное число нулей. Более того, существуют такой номер $n_0$ и число $a_0 > 0$, что при всех $n \geq n_0, \; \Lambda_n \geq a_o n^{-2}$.}

Наряду с точкой $t_0$ введем в рассмотрение еще функцию $\psi(x)$, удовлетворяющую условиям \eqref{eq2.1} (выше мы показали, что такое предположение допустимо). И пусть ее разложение в ряд Фурье имеет вид:
\begin{equation}
\psi(x) = \sum_{k=0}^{\infty}\psi_k \sin{kx}
\label{eq2.3}
\end{equation}
Задача 2 состоит в нахождении функции $f \in C^1 ([0,\pi]), f(0) = f(\pi) = 0$, при которой решение $u_{\omega} (x,t)$ задачи \eqref{s1}-\eqref{s2} удовлетворяет соотношению:
\begin{equation}
\bigl\| u_{\omega}(x, t_0) - \psi(x) \bigr\|_{C([0,\pi])} = o(1), \omega \to \infty.
\label{eq2.4}
\end{equation}

Для формулировки основного результата данного пункта обозначим номера нулевых значений $\Lambda_n$ через $n_1, n_2 , ... ,n_s$, а через $M_0$ -- множество этих номеров.

\begin{theorem} \qquad

1. Задача 2 имеет единственное решение тогда и только тогда, когда множество $M_0$ пусто. При этом $f = \sum_{n=1}^{\infty}f_n\sin{nx}$, где $f_n = \frac{\psi_n}{\Lambda_n}$.

2. Если множество $M_0$ не пусто, то задача 2 разрешима в том и только том случае, когда $\psi_{n_j} = 0, n_j \in M_0.$ При этом коэффициенты $f_n = \frac{\psi_n}{\Lambda_n}, n \not\in M_0,$ а $f_n$ при $n \in M_0$ - произвольные числа.

\end{theorem}

Укажем простые естественные достаточные ограничения на $t_0$, обеспечивающие единственность $f$. Если $r_0(0) \neq 0$, то в качестве $t_0$ можно взять любое число на участке (0,1], не превышающее первый нуль функции $r_0(t), t \in [0,1]$. Если же $r_0(0) = 0$, то в качестве $t_0\in (0,1)$ можно взять такое число, что $r_0(t_0) \neq 0$ и $r_0(t)$ на участке $t \in [0, t_0]$ не меняет знака.

\noindent\textbf{$2.3^{\circ}.$ Доказательство основных результатов параграфа}

\noindent\textbf{Доказательство теоремы 3} \newline
Рассмотрим функцию
\begin{center}
$W_{\omega}(x,t) = u_{\omega}(x,t) - u_0(x,t) = \sum_{n=1}^{\infty}f_n\sin{nx} \int\limits_0^t e^{-n^2(t-s)} r_1 (s, \omega s) ds = \sum_{n=1}^{\infty}f_n\sin{nx} \int\limits_0^t e^{-n^2(t-s)} r_1 (s, \omega s) ds + \sum_{n=n_0 + 1}^{n_0}f_n\sin{nx} \int\limits_0^t e^{-n^2(t-s)} r_1 (s, \omega s) ds \equiv S_{\omega, 1} + S_{\omega, 2}, n_o \in \mathbb{N}$
\end{center}
Пусть $\varepsilon$ - произвольное положительное число. Учитывая, что $\int\limits_0^t e^{-n^2 (t-s)} r_1 (s,\omega s) ds = O (n^{-2}), n \to \infty$, равномерно относительно $t \in [0,1], \omega > 0$, подберем $n_0$ настолько большим, что при всех $\omega > 0$
\begin{center}
$\bigl\| S_{\omega, 2} \bigr\|_{C(\overline{\Pi})} < \frac{\varepsilon}{2}.$
\end{center}
После этого, в силу равенства $\left\langle r_1 (t,\tau) \right\rangle_{\tau} = 0$, подберем $\omega_0$ столь большим, что при найденном $n_0$ и всех $\omega > \omega_0$
\begin{center}
$\bigl\| \sum_{n=1}^{n_0}f_n\sin{nx} \int\limits_0^t e^{-n^2(t-s)} r_1 (s, \omega s) ds \bigr\|_{C(\overline{\Pi})} < \frac{\varepsilon}{2}.$
\end{center}
Таким образом, при $\omega > \omega_0$
\begin{center}
$\bigl\| W_{\omega} \bigr\|_{C(\overline{\Pi})} < \varepsilon.$
\end{center}
Теорема доказана.

\noindent\textbf{Доказательство теоремы 4.} \newline
Предположим, что функция $f \in C^1([0, \pi])$ найдена. Тогда в силу теоремы 3 и условий \eqref{eq2.3}, \eqref{eq2.4} справедливо равенство
\begin{center}
$\sum_{n=1}^{\infty}f_n\Lambda_n\sin{nx} = \sum_{n=1}^{\infty}\psi_n\sin{nx}$
\end{center}
Ввиду предположений \eqref{eq2.1} найдутся константы $b_n > 0, n = 1,2,\ldots$, такие что при всех натуральных $n$
\begin{center}
$\left|\psi_n\right| \leq b_n n^{-4}, \sum_{n = 1}^{\infty} |b_n|^2 < \infty,$
\end{center}
Отсюда в силу леммы при $n \geq n_0$
\begin{center}
$\left|f_n\right| \leq c_n n^{-2}, \; c_n = b_n a_0,$
\end{center}
так что $f \in C^1 ([0,\pi])$. Если все $\Lambda_n \neq 0$, то однозначно находим
\begin{center}
$f = \sum_{n=1}^{\infty}f_n\sin{nx} \in C^1([0,\pi]), \; f_n = \frac{\psi_n}{\Lambda_n}.$
\end{center}
Остальная часть теоремы 4 очевидна.

\section{Известно лишь среднее значение зависящего от времени сомножителя источника}
\label{sect3}
\subsection*{$3.1^{\circ}$ Задача 3}

\qquad В этом разделе вновь рассмотрим систему вида \eqref{s1}-\eqref{s2} с $f(x,t) \equiv f(x), f \in C^3 ([0,\pi]), f^{(2k)} (0) = f^{(2k)} (\pi) = 0, k = 0,1$, и функцией $r = r_0 + r_1$, принадлежащей классу (А). Будем здесь считать, что функция $r_0$ известна, а функции $f$ и $r_1$ не известны. Пусть точка $t_0 \in (0,T]$, а $\Lambda_n = \Lambda_n (t_0), n=1,2,\ldots$ -- последовательность чисел, вычисленных в соответствие с формулой \eqref{eq2.11}. Для простоты будем предполагать, что $\Lambda_n \neq 0, n = 1,2,\ldots$ Достаточные условия для этого указаны в п. $2.2^{\circ}$. Зададим четыре функции: $\psi(x), \varphi_0 (t), \varphi_1 (t)$ и $\varphi_2 (t,\tau)$ следующим образом. Пусть функция $\psi \in C^{6} ([0,\pi]), \psi^{(2n)} (0) = \psi^{(2n)} (\pi) = 0, n = 0,1,2.$ Положим
\begin{equation}
\widetilde{f} (x) = \sum_{n = 1}^{\infty} \widetilde{f}_n \sin{nx}, \; \widetilde{f}_n = \frac{\psi_n}{\Lambda_n}
\label{fw}
\end{equation}
Тогда $\widetilde{f} \in C^3([0,\pi]), \widetilde{f}^{(2k)} (0) = \widetilde{f}^{(2k)} (\pi) = 0, k = 0,1.$ Пусть $x_0 \in (0,\pi)$ -- точка, в которой $\widetilde{f}(x_0) \neq 0$. Пусть функция $\varphi_0 (t)$ задана равенствами
\begin{equation}
\widetilde{f}(x_0) r_0 (t) + \int\limits_0^t K(t-s) r_0 (s) ds = \varphi'_0 (t), \varphi_0 (0) = 0,
\label{eq3.2}
\end{equation}
где $K(t,s)$ -- та же функция, что в $\S 1$, а функции $\varphi_1, \varphi_2$ удовлетворяют условиям п.$1.2^{\circ}$ с заменой (в выражении $\varphi_1$ через $\varphi_2$) функции $f(x)$ на $\widetilde{f}(x)$. Эти условия и \eqref{eq3.2} можно назвать условиями согласования функций $r_0, \psi, \varphi_0, \varphi_1$ и $\varphi_2$.

Из результатов разделов 1, 2 вытекает следующая теорема.

\begin{theorem}
Пусть функции $r_0, \psi, \varphi_0, \varphi_1,$ и $\varphi_2$ и точки $x_0, t_0$ удовлетворяют указанным в этом пункте условиям. Тогда существует единственная (в указанных в п.$3.1^{\circ}$ классах) пара функций $f$ и $r_1$, при которых решение $u_{\omega} (x,t)$ задачи \eqref{s1}-\eqref{s2} при $f(x,t) \equiv f(x)$ удовлетворяет условиям \eqref{babich-inv} и \eqref{eq2.4}. При этом функция $f(x) = \widetilde{f}(x)$ вычисляется по формуле \eqref{fw}, а $r_1 (t,\tau) = (f(x_0))^{-1} \frac{\partial}{\partial \tau} \varphi_2 (t,\tau)$. Решение задачи 3 сводится к решению уравнения Вольтерры второго рода.
\end{theorem}

\section{Не известны оба сомножителя источника}
\label{sect4}
\subsection*{$4.1^{\circ}$ Задача 4}

\qquad Пусть $\Pi$ и $\Gamma$ -- те же, что в $\S1$. Рассмотрим систему \eqref{s1}-\eqref{s2}, в которой функции $f(x,t) \equiv f(x)$ и $r(t,\tau)$ не заданы. Однако известно, что $r(t,\tau)$ -- функция класса (А), а $f(x,t) \equiv f(x) = \sum\limits_{n=1}^N f_n \sin nx$ -- функция с заданным числом $N$ гармоник с неизвестными амплитудами $f_n$. Пусть, кроме того, заданы точки $t_0 \in (0,T), x_j \in (0,\pi), j = \overline{0,N-1}, x_i \neq x_k$ при $i \neq k$, а также функции $\varphi_0(t), \varphi_2 (t,\tau), \varphi_0 \in C^1 ([0,T]), \varphi_0(0) = 0, \varphi_2 \in C^{1+\alpha, 0} (S), \alpha \in (0,1), \varphi_2(t,\tau)$ -- $2\pi$-периодична по $\tau$ и функции $\alpha_j \in C^1 ([t_0-\delta,t_0+\delta]), \delta >0, j = \overline{1,N-1}, (t_0-\delta,t_0+\delta) \subset (0,T)$. Символом $\alpha (t)$ будем обозначать вектор-функцию с компонентами $\varphi_0 (t), \alpha_j (t), t \in (t_0 - \delta, t_0 + \delta), j=\overline{0,N-1}$. Аналогичные обозначения будем использовать и в случае других вектор-функций и векторов.

Введем еще функцию, точнее, представление:
\begin{equation*}
\varphi_1 (t) = -\frac{1}{f(x_0)} \varphi_2 (0,0) \sum\limits_{n=1}^N f_n \sin nx_0 e^{-n^2 t}, \, t \in [0,T], \, f(x_0) \neq 0,
\end{equation*}
в котором фигурирует пока неизвестная функция $f(x)$.

Задача 4 состоит в определении таких функций $r(t,\tau)$ и $f(x)$, удовлетворяющих указанным выше условиям, при которых для решения $u_{\omega} (x,t)$ задачи \eqref{s1}, \eqref{s2} выполнены асимптотические формулы
\begin{equation}
\left\| u_{\omega} (x_0,t) - [\varphi_0 (t) + \omega^{-1}( \varphi_1 (t) + \varphi_2 (t,\omega t) )]  \right\|_{C([0,T])} = o(\omega^{-1}),
\end{equation}
\begin{equation}
\left\| u_{\omega} (x_j,t) - \alpha_j (t) \right\|_{C([t_0-\delta,t_0+\delta])} = o(1), j=\overline{1,N-1}, \omega \to \infty.
\end{equation}

Прежде чем сформулировать основной результат данного параграфа введем ряд обозначений и сформулируем некоторые дополнительные условия. Рассмотрим систему уравнений
\begin{equation}
\left\{ \begin{array}{c}
\sum\limits_{n=1}^N \psi_n \sin nx_0 = \varphi_0 (t_0), \\
\sum\limits_{n=1}^N \psi_n \sin nx_j = \alpha_j (t_0), j = \overline{1,N},
\end{array} \right.
\label{s4.3}
\end{equation}
относительно неизвестных $\psi_n, n=\overline{1,N}$. Поскольку матрица $A = (\sin nx_j)_{n = 1, j = 0}^{N, N-1}$ невырождена, то из \eqref{s4.3} однозначно найдем вектор $\psi \equiv \psi (\alpha(t_0))$. Рассмотрим затем систему
\begin{equation}
\left\{
\begin{array}{c}
\sum\limits_{n=1}^N \varphi_n \sin nx_0 = \sum\limits_{n=1}^N n^2 \psi_n \sin nx_0 + \varphi'_0 (t_0), \\
\sum\limits_{n=1}^N f_n \sin nx_j = \sum\limits_{n=1}^N n^2 \psi_n \sin nx_j + \alpha'_j (t_0), j = \overline{1,N-1},
\end{array}
\right.
\label{s4.4}
\end{equation}
относительно неизвестных $f_n, n=\overline{1,N},$ из которой однозначно найдем вектор
\begin{equation}
f \equiv F_1 (\psi (\alpha(t_0)), \alpha'(t_0)) \equiv F (\alpha, t_0) \equiv F.
\end{equation}
Будем предполагать, что выполнено соотношение
\begin{equation}
\sum\limits_{n=1}^N f_n \sin nx_0 \neq 0.
\label{eq4.6}
\end{equation}
Нормируем искомую функцию $r_0 (t)$ условием $r_0 (t_0) = 1$. Рассмотрим теперь уравнение Вольтерры второго рода
\begin{equation}
f(x_0) l(t) + \int\limits_0^t K(t,s) l(s) ds = \mu (t), \mu \in C([0,T]),
\label{eq4.7}
\end{equation}
где $K(t,s) = - \sum\limits_{n=1}^N n^2 F_n e^{-n^2(t-s)} \sin nx_0, \, f(x) = \sum\limits_{n=1}^N f_n \sin nx$. Его единственное в пространстве $C([0,T])$ решение $l(t)$ обозначим символом $S (\alpha, t_0, x_0, \mu (t))$.

\begin{theorem}
Для любого набора функций $\varphi_i, i=\overline{0,2}, \alpha_j, j=\overline{1,N-1}$, и точек $t_0,x_k, k=\overline{0,N-1}$, удовлетворяющих указанным выше условиям, задача 4 однозначно разрешима тогда и только тогда, когда выполнены следующие условия согласования:
\begin{equation}
\sum\limits_{n=1}^N \sin nx_j F_n (\alpha, t_0) \int\limits_0^t e^{-n^2 (t-s)}  S (\alpha, t_0, x_0, \varphi'_0(s)) ds = \alpha_j (t), j=\overline{1,N-1}, t \in (t_0 - \delta, t_0 + \delta).
\label{eq4.8}
\end{equation}
При их выполнении для нахождения функции $f(x)$ требуется решить две системы линейных уравнений с единой невырожденной основной матрицей, для нахождения $r_0 (t)$ -- уравнения Вольтерры второго рода \eqref{eq4.7} с $\mu(t) = \varphi'_0(t)$, а $r_1 (t,\tau) = (f(x_0))^{-1} \frac{\partial}{\partial\tau} \varphi_2 (t,\tau)$.
\end{theorem}

Доказательство теоремы 6 мы опускаем; оно основано на тех же соображениях, что использованы выше. Из того факта, что условия согласования \eqref{eq4.8} является и необходимым для разрешимости задачи 4, следует, что это равенство выполняется на непустом (даже широком) множестве исходных данных. Тем не менее приведем один простой иллюстративный пример к теореме 6.

\subsection*{$4.2^{\circ}$ Пример}

\qquad Рассмотрим задачу \eqref{s1}-\eqref{s2} в условиях теоремы 6 с $N=2$. В качестве исходных данных возьмем следующие:
\begin{equation}
\begin{array}{c}
t_0 = 1, x_0 = \frac{\pi}{2}, x_1 = \frac{\pi}{6},\\
\varphi_0 (t) = e^{-t} + t - 1, \\
\alpha_1 (t) = \frac{1}{2} (t + e^{-t} - 1) + \frac{\sqrt{3}}{32} (4t + e^{-4t} - 1), \\
\varphi_2 (t,\tau) = - \cos{\tau}.
\end{array}
\end{equation}

Проверим для этих данных справедливость условия согласования \eqref{eq4.8}. Система \eqref{s4.3} в данном случае принимает вид
\begin{equation*}
\left\{
\begin{array}{c}
\psi_1 = e^{-1}, \\
\frac{1}{2} \psi_1 + \frac{\sqrt{3}}{2}\psi_2 = \frac{1}{2 e} + \frac{\sqrt{3}}{32} (3+ e^{-4}),
\end{array}
\right.
\label{examp1}
\end{equation*}
откуда находим $\psi_1 = e^{-1}, \psi_2 = \frac{1}{16} (3 + e^{-4})$. Система \eqref{s4.4} имеет вид
\begin{equation*}
\left\{
\begin{array}{c}
f_1 = - e^{-1} + 1 + e^{-1}, \\
\frac{1}{2} f_1 + \frac{\sqrt{3}}{2}f_2 = \frac{1}{2} + \frac{\sqrt{3}}{2}
\end{array}
\right.
\label{examp2}
\end{equation*}
откуда определим $f_1 = 1, f_2 = 1$. Справедливо соотношение \eqref{eq4.6}: $\sin \frac{\pi}{2} - \sin \pi = 1 \neq 0$. Из уравнения Вольтерры второго рода
\begin{equation*}
r_0 (t) - \int\limits_0^t e^{-(t-s)} r_0 (s) ds = - e^{-t} + 1,
\label{examp3}
\end{equation*}
и находим $r_0 (t) = t$. Теперь можно выписать условие согласования \eqref{eq4.8}:
\begin{equation*}
\frac{e^{-t}}{2} \int\limits_0^t e^s s ds + \sqrt{3}\frac{e^{-4t}}{2} \int\limits_0^t e^{4s} s ds = \alpha_1 (t), t \in (t_0 - \delta, t_0 + \delta), \delta > 0.
\label{examp4}
\end{equation*}
Легко проверяется его справедливость. Итак, имеют место представления:
\begin{equation*}
f(x) = \sin x - \sin 2x, \, r_0 (t) = t, \, r_1 (t,\tau) = \sin \tau.
\label{examp5}
\end{equation*}

\section{Заключение}
\label{sect5}

Решены четыре задачи о восстановлении высокочастотного источника в одномерном уравнении теплопроводности с однородными начально-краевыми условиями по некоторым сведениям о частичной асимптотике его решения. Показано, что источник удается полностью восстановить по определенным данным о неполной (двучленной) асимптотике решения. Постановке каждой задачи о восстановлении источника в работе предшествует построение с обоснованием асимптотики решения исходной начально-краевой задачи.


\newpage

\begin{thebibliography}{text}

\bibitem{babich-bib1}
{\it Денисов А.\,М.\/} Асимптотика решений обратных задач для гиперболического уравнения с малым параметром при старшей производной // ЖВМиМФ. 2013. Т.~53, 5. С.~744--752.

\bibitem{babich-bib2}
{\it Денисов А.\,М.\/} Задачи определения неизвестного источника в параболическом и гиперболическом уравнении // ЖВМиМФ. 2015. Т.~55, 4. С.~830--835.

\bibitem{babich-bib3}
{\it Babich P.\,V., Levenshtam V.\,B.\/} Direct and inverse asymptotic problems high-frequency terms // Asymptotic Analysis. 2016. Т.~97. С.~329--336.

\bibitem{babich-bib31}
{\it Зеньковская С.\,М., Симоненко И.\,Б./} О влиянии вибрации высокой частоты на возникновение конвекции // Механика жидкости и газа. 1966. 5. C.~51--55.

\bibitem{babich-bib32}
{\it Симоненко И.\,Б.\/} Обоснование метода усреднения для задачи конвекции в поле быстро осциллирующих сил и для других параболических уравнений // Мат. сб. 1972. Т.~87(129), 2. C.~236--253.

\bibitem{babich-bib33}
{\it Левенштам В.\,Б.\/} Метод усреднения в задаче конвекции при высокочастотных наклонных вибрациях // Сиб. матем. журн. 1996. Т.~37, 5. C.~1103--1116.

\bibitem{babich-bib5}
{\it Левенштам В.\,Б.} Асимптотическое разложение решения задачи о вибрационной конвекции // ЖВМиМФ. 2000. Т.~40, 9. С.~1416--1424.

\bibitem{babich-bib4}
{\it Вишик М.\,И., Люстерник Л.\,А.} Регулярное вырождение и пограничный слой для линейных дифференциальных уравнений с малым параметром // УМН. 1957. Т.~12, 5. С.~3--122.

\bibitem{babich-bib6}
{\it Левенштам В.\,Б.} Асимптотическое интегрирование параболических задач с большими высокочастотными слагаемыми // Сиб. матем. журн. 2005. Т.~46, 4. С.~805--821.

\bibitem{babich-bib7}
{\it Левенштам В.\,Б.} Построение старших приближений метода усреднения для параболических начально-краевых задач методом пограничного слоя // Изв. вузов. Матем. 2004. 3. С.~41--45.

\end{thebibliography}
\end{document}